# DISCUSSION OF "BREAKDOWN AND GROUPS" BY P. L. DAVIES AND U. GATHER

By Peter J. Rousseeuw

*Universiteit Antwerpen*

**1. General comments.** The interesting paper of Davies and Gather (henceforth [DG]) pulls together results on upper bounds on the breakdown value of translation equivariant location estimators [Donoho (1982)], regression estimators [Rousseeuw (1984)] and affine equivariant scatter estimators [Davies (1987)] into a single framework of group equivariance. I can only agree with them on the important role of the latter notion in obtaining nontrivial bounds. [By the way, I prefer the term breakdown *value* myself because it is not a point, and the term "value" captures both its dimension (one) and its orientation (we aim for higher, not lower values).]

The theory in [DG] is formulated for estimators that are uniquely defined, but it seems to work just as well in the general case. Then $T(P)$ is a set, and we can follow the implicit convention of saying that it breaks down when any member of $T(P)$ does. We only need to redefine $D(T(P), T(Q))$ in (2.4) as a supremum over all pairs of members of $T(P)$ and $T(Q)$.

The new applications of the theory are fascinating, for example, to the Michaelis–Menten model (with a nontrivial bound) and logistic regression (without one). I am less convinced by the fragility argument illustrated by the difference between the contaminated samples (6.2) and (6.3). It is true that this proof of the upper bound only covers (6.2), but in some sense that is enough since breakdown is a worst-case concept and the bound is not specific for the median but for all translation estimators. But anyway, the behavior of the median at (6.3) can be derived from that at (6.2) by a variety of other properties that it possesses. For instance, its monotonicity property alone suffices. Or we can use the property that when you move the observations over distances of at most $\delta$, the median changes by at most $\delta$. This holds for any $\delta > 0$, and is a Lipschitz property for the metric on samples defined as

$$(1.1) \qquad d((x_1, \ldots, x_n), (y_1, \ldots, y_n)) = \min_{\pi \in S_n} \max_i |x_i - y_{\pi(i)}|,$$

where $S_n$ is the set of all permutations on $\{1,\ldots,n\}$. Note that (1.1) is equal to $\max_i |x_{i:n} - y_{i:n}|$ [see Rousseeuw and Leroy (1987), pages 127–128]. People who compute maxbias curves always use the properties of the actual estimator. Perhaps we should not expect much more elegant results for contaminated samples that break down the estimator than for those that yield a finite bias.

**2. The maximal breakdown value of affine equivariant location estimators.** From here on I will focus on the open problem in Section 5.2 of [DG]. It has been known since Donoho (1982) that the finite-sample breakdown value (fsbv) of translation equivariant estimators of location is at most $\lfloor (n+1)/2 \rfloor /n$ and that this bound is sharp. The bound obviously holds also for affine equivariant location estimators, but it may not be sharp for them. In one dimension ($k=1$) it is, but for $k \geq 2$ this has been an open problem for over 20 years. During that time many affine location estimators were constructed with an fsbv of $\lfloor (n-k+1)/2 \rfloor /n$, such as the MVE and MCD of Rousseeuw (1984), location $S$-estimators [Davies (1987), Rousseeuw and Leroy (1987)] and a modification of the Stahel–Donoho estimator [Tyler (1994), Gather and Hilker (1997)]. Since $\lfloor (n-k+1)/2 \rfloor /n$ is known to be the sharp upper bound for affine equivariant scatter estimators [Davies (1987)], it has seemed plausible that it could also be the upper bound for affine location. Over the years there have been several attempts to attain the upper bound $\lfloor (n+1)/2 \rfloor /n$, but as far as I know none has succeeded. [The result in Zuo (2004) does not count because it uses a weaker version of the fsbv which requires that all the contaminating points coincide.]

Let us consider any data set $X = \{x_1,\ldots,x_n\} \subset \mathbf{R}^k$ (from here on always $k \geq 2$ and $n > k$) which is in general position (GP). By GP we mean that no more than $k$ data points lie on any affine hyperplane. This holds a.s. when sampling from an absolutely continuous distribution. The convex hull $\text{conv}(X)$ is then a polytope with faces that contain exactly $k$ data points. [In $\mathbf{R}^3$ the faces are two-dimensional, and in general they are $(k-1)$-dimensional.] Note that $\text{conv}(X)$ can be stretched arbitrarily by replacing even a single point of $X$ by an outlier. Since we are studying very robust estimators $T$, it is natural to require that $T$ should not lie on the boundary of $\text{conv}(X)$ or outside of $\text{conv}(X)$. A slightly more general formulation of this requirement is the following condition:

($C_h$) Let $X = \{x_1,\ldots,x_n\} \subset \mathbf{R}^k$ be in general position, with $n > k \geq 2$. Let $u$ be a direction such that the inner products $y_i = u'x_i$ satisfy $y_1 = \cdots = y_h < y_{h+1} \leq \cdots \leq y_n$ (after renumbering) for the specified number $h$, with $1 \leq h \leq k$. Then there exists an $\alpha > 0$ (which depends only on $k$ and the $y_1,\ldots,y_n$) such that $u'T(X) \geq y_h + \alpha$.

The typical case is to take $h = k$. For any face of $\text{conv}(X)$ we can take the orthogonal direction $u$ pointing to the inside of $\text{conv}(X)$, so Condition ($C_k$) says that $T$ cannot lie on or arbitrarily close to the boundary of $\text{conv}(X)$ or outside of it. [Note that $\text{conv}(X)$ is the intersection of halfspaces containing $X$ and having a face of $\text{conv}(X)$ on their boundary.] For $h < k$ the condition becomes somewhat weaker; for example, Condition ($C_1$) only says that $T$ cannot come arbitrarily close to a vertex of $\text{conv}(X)$ or lie outside of $\text{conv}(X)$.

Condition ($C_k$) is intuitive for a robust estimator. For instance, Condition ($C_k$) holds for all estimators that can be written as a weighted mean $(\sum_i w_i x_i)/(\sum_i w_i)$ where $0 \leq w_i \leq 1$ and at least $k+1$ of the $w_i$ equal 1 [it suffices to put $\alpha = (y_{k+1} - y_k)/n$]. This encompasses, for example, trimmed means and the minimum covariance determinant estimator (MCD). Moreover, a robust estimator would typically be expected to have a reasonably large Tukey depth, for example,

$$\text{depth}(T, X) \geq k + 1 \tag{2.2}$$

(at least for large enough $n$, when there are many depth contours). Condition (2.2) implies Condition ($C_k$) and is another way of saying that $T$ should not be in the outskirts of the data cloud.

THEOREM 1. *Consider a data set $X = \{x_1, \ldots, x_n\} \subset \mathbf{R}^k$ in general position with $n > k$. Let $T$ be an affine equivariant location estimator satisfying Condition ($C_h$) with $1 \leq h \leq k$. Then $\text{fsbv}(T, X) \leq \lfloor (n - h + 1)/2 \rfloor / n$.*

PROOF. Put $\hat{\theta} := T(X) \in \mathbf{R}^k$. Since $X$ is in GP, $\text{conv}(X \cup \{\hat{\theta}\})$ has at least one face not containing $\hat{\theta}$. Take an $h$-subset $S$ of the $k$ data points on this face. Then there exists an affine hyperplane $L$ which contains $S$ such that both $\hat{\theta}$ and $X \setminus S$ lie strictly on the same side of $L$. Assume w.l.o.g. that $0 \in L$. Denote the unit normal vector to $L$ in the direction of $X \setminus S$ as $e_1$ and take an orthonormal basis $\{e_2, \ldots, e_k\}$ of $L$. After renumbering, the $x_{i1} := e_1' x_i$ satisfy $0 = x_{1,1} = \cdots = x_{h,1} < x_{h+1,1} \leq \cdots \leq x_{n,1}$, hence $\hat{\theta}_1 \geq \alpha > 0$ by Condition ($C_h$).

Let us assume that $T$ cannot be broken down by replacing any $m$-subset $B$ of $X$, where $m = \lfloor (n - h + 1)/2 \rfloor$, by an arbitrary $m$-set $B'$ yielding the contaminated data set $X' := (X \setminus B) \cup B'$. This means that there exists a finite radius $M$ such that for any contaminated data set $X'$ of this type it holds that $T(X') \in B(\hat{\theta}, M) := \{x \in \mathbf{R}^k; \|x - \hat{\theta}\| \leq M\}$.

We will now construct a linear transformation which leaves $S$ invariant and moves $X \setminus S$ as well as $\hat{\theta}$. For this we consider the "shear transform" $g_\gamma$ given by the nonsingular matrix

$$\begin{bmatrix} 1 & 0 & 0 \\ \gamma & 1 & 0 \\ 0 & 0 & I_{k-2} \end{bmatrix} \tag{2.3}$$

relative to the basis $\{e_1, \ldots, e_k\}$, for $\gamma \in \mathbf{R}$. We note that $g_\gamma(e_j) = e_j$ for all $j \neq 1$, hence $g_\gamma(x_i) = x_i$ for $i = 1, \ldots, h$, but at the same time $g_\gamma(e_1) = e_1 + \gamma e_2$. Denoting $\hat{\theta} = (\hat{\theta}_1, \ldots, \hat{\theta}_k)^T$, we find $g_\gamma(\hat{\theta}) = (\hat{\theta}_1, \hat{\theta}_2 + \gamma\hat{\theta}_1, \hat{\theta}_3, \ldots, \hat{\theta}_k)^T$ with $\hat{\theta}_1 > 0$, hence $\|g_\gamma(\hat{\theta}) - \hat{\theta}\| = |\gamma|\hat{\theta}_1$ goes to infinity for increasing $\gamma$. Analogously, the image of any data point $x_i$ with $i = h+1, \ldots, n$ is of the form $g_\gamma(x_i) = x_i + \gamma x_{i1} e_2$, so all $g_\gamma(x_i)$ move in the direction of $e_2$ and $(g_\gamma(x_i))_1 = x_{i1}$. Each point travels a distance $\|g_\gamma(x_i) - x_i\| = |\gamma||x_{i1}| \geq |\gamma||x_{h+1,1}|$.

Let us partition $X \setminus S$ into two sets $A$ and $B$ with $|B| = m = \lfloor (n - h + 1)/2 \rfloor$ and $|A| = n - h - |B|$. (If $n - h$ is even, we find $|A| = |B|$, whereas for odd $n - h$ we have $|A| = |B| - 1$.) We will replace $B$ by $B_\gamma := g_\gamma(B)$ yielding the contaminated data set $X'_\gamma := S \cup A \cup B_\gamma$. Note that $X'_\gamma$ is in GP for all but a finite number of $\gamma$ values. Put $\Gamma = \{\gamma; X'_\gamma \text{ is in GP}\}$. For all $\gamma \in \Gamma$ it holds that $T(X'_\gamma) \in H := \{z \in \mathbf{R}^k; z_1 \geq \alpha\}$ by Condition ($C_h$).

For any $\gamma$ the image of $B(\hat{\theta}, M)$ through $g_\gamma$ is an ellipsoid with center $g_\gamma(\theta)$. For a large enough $\gamma \in \Gamma$ it holds that $B(\hat{\theta}, M) \cap g_\gamma(B(\hat{\theta}, M)) \cap H = \varnothing$. We know that $T(X'_\gamma) \in B(\hat{\theta}, M)$ by assumption. On the other hand, we can also write $X'_\gamma = g_\gamma(S \cup A_{-\gamma} \cup B)$, which implies $T(X'_\gamma) \in g_\gamma(B(\hat{\theta}, M))$. Since $T(X'_\gamma) \in H$ it follows that $T(X'_\gamma) \in B(\hat{\theta}, M) \cap g_\gamma(B(\hat{\theta}, M)) \cap H = \varnothing$. This contradiction proves the desired upper bound on fsbv. $\square$

In the typical case where $h = k$, Theorem 1 yields the upper bound $\lfloor (n - k + 1)/2 \rfloor/n$ which has been attained. This says that any affine location estimator $T$ with a higher fsbv must be somewhat strange in the sense of not satisfying Condition ($C_k$), so $T$ can be arbitrarily close to the boundary of $\text{conv}(X)$ or even lie outside it. Any $T$ which were to attain the translation equivariant bound $\lfloor (n+1)/2 \rfloor/n$ cannot even satisfy Condition ($C_1$), so at times it must be arbitrarily close to a vertex of $\text{conv}(X)$ or lie outside it. It is counterintuitive that an estimator with maximal fsbv would have such a low Tukey depth (at most 1).

So far the only published result with higher fsbv than $\lfloor (n - k + 1)/2 \rfloor/n$ is the projection median (PM) of Zuo (2003), which attains $\lfloor (n - k + 2)/2 \rfloor/n$ by using a univariate scale estimator $\text{MAD}_{k-1}$ in its definition. By Theorem 1, this estimator cannot satisfy Condition ($C_k$). Here is a bivariate counterexample (which can be extended to $\mathbf{R}^k$). Start with the data points $z_1 = (0, \delta)$ and $z_2 = (0, -\delta)$ for some $\delta > 0$. Add $m$ points $(x_i, y_i)$ with $x_i$ equispaced between 10 and 20 and $y_i = x_i + \delta u_i$ where the noise $u_i$ is such that these points are in GP. Add another $m$ points with the same $x_i$ but with $-y_i$. Then the $n = 2m + 2$ points of $Z$ are in GP for all but finitely many $\delta$. When $\delta \to 0$, the outlyingness $\text{Out}(0,0)$ tends to the outlyingness of 0 relative to $\{0, 0, x_1, x_1, \ldots, x_m, x_m\}$; hence for any $0 < \delta < 1$ we have $\text{Out}(0, 0) < M$ for some $M < \infty$. We will prove that for any $\varepsilon > 0$ there is a $\delta_0 > 0$ such that $\delta < \delta_0$ implies $\|\text{PM}(Z)\| < \varepsilon$. By projecting in the direction orthogonal to $y = -x$ we see that $\text{MAD}_1$ tends to 0, so for small enough $\delta$ all points (not necessarily data points) in $\mathbf{R}^2$ lying farther than $\varepsilon/\sqrt{2}$ away from the line $y = -x$ have $\text{Out} > M$. The same holds for points farther than $\varepsilon/\sqrt{2}$ from $y = x$. Therefore $\text{PM} \to (0,0)$; hence $\alpha$ in Condition $(\text{C}_k)$ is zero.

Note that Theorem 1 fits in the framework of [DG] with $G$ the affine group on $\mathbf{R}^k$. The main difference is that here we first fix a set $B$ (our $h$-subset) and then a subgroup of $G$ which keeps $B$ invariant, whereas condition (3.3) in [DG] is over many possible $B$. Afterward we put $g := g_1$ [i.e., (2.3) with $\gamma = 1$], yielding $\Delta(P_n) = h/n$. The remainder of the proof of Theorem 3.1 in [DG] can then be retraced by noting that for any integer $m$ it holds that $g^m = g_m$ (the shear transform with $\gamma = m$). We basically set aside $h$ points and then apply our usual reasoning to the remaining $n - h$ points.

Also note that Condition $(\text{C}_h)$ and Theorem 1 can be extended to situations without general position. As long as $T$ satisfies Condition $(\text{C}_h)$ without the GP condition (this is a stronger assumption), and $X$ does have $h$ points whose inner products with some $u$ satisfy $y_1 = \cdots = y_h < y_{h+1} \leq \cdots \leq y_n$, the upper bound $\text{fsbv}(T, X) \leq \lfloor (n - h + 1)/2 \rfloor / n$ holds. In this situation it is even allowed that $h > k$ (which could not happen under GP).

**Acknowledgment.** I would like to thank Yijun Zuo for stimulating discussions.

Department of Mathematics
  and Computer Science
Universiteit Antwerpen
Middelheimlaan 1
B-2020 Antwerp
Belgium
e-mail: Peter.Rousseeuw@ua.ac.be